\documentclass[11pt]{amsart}

\usepackage[colorlinks = true,
            linkcolor = red,
            urlcolor  = magenta,
            citecolor = black,
            anchorcolor = blue]{hyperref}
\usepackage{color}
\usepackage[shortalphabetic,initials,nobysame]{amsrefs}
\usepackage{upgreek}
\usepackage{tikz}
\usepackage[applemac]{inputenc} 
\usetikzlibrary{arrows}
\usepackage{xcolor}
\usepackage[all]{xy}
\usepackage{graphicx}
\usepackage{wrapfig}
\usepackage{yfonts}
\usepackage{graphicx}
\usepackage{amssymb}
\usepackage{epstopdf}
\usepackage{mathrsfs}
\usepackage[mathcal]{eucal}
\usepackage{bm}

\renewcommand{\eprint}[1]{\href{https://arxiv.org/abs/#1}{#1}}

\DeclareMathOperator{\Lie}{Lie}

\BibSpec{article}{%
    +{}  {\PrintAuthors}                {author}
    +{,} { \textit}                     {title}
     +{,} {}                            {note}
    +{.} { }                            {part}
    +{:} { \textit}                     {subtitle}
    +{,} { \PrintContributions}         {contribution}
    +{.} { \PrintPartials}              {partial}
    +{,} { }                            {journal}
    +{}  { \textbf}                     {volume}
    +{}  { \PrintDate}                {date}
    +{,} { \issuetext}                  {number}
    +{,} { \eprintpages}                {pages}
    +{,} { }                            {status}
    +{,} { \DOI}                   {doi}
    +{,} { \eprint}        {eprint}
      +{,} {\publisher}                {publisher}
    +{,} { \address}                   {address}
    +{}  { \parenthesize}               {language}
    +{}  { \PrintTranslation}           {translation}
    +{;} { \PrintReprint}               {reprint}
    +{.} {}                             {transition}
    +{}  {\SentenceSpace \PrintReviews} {review}
}

\newtheorem{Thm}{Theorem}[section]

\newtheorem{Prop}[Thm]{Proposition}
\newtheorem{Cor}[Thm]{Corollary}
\newtheorem{Con}[Thm]{Conjecture}

\theoremstyle{definition}
\newtheorem{Def}[Thm]{Definition}
\theoremstyle{remark}
\newtheorem{Rem}[Thm]{Remark}

\newtheoremstyle{named}{}{}{\itshape}{}{\bfseries}{.}{.5em}{#1 #3}
\theoremstyle{named}

\def\si{{\Sigma}}

\def\fb{\mathfrak{b}}
\def\g{\mathfrak{g}}
\def\Frenkel:2013uda{\mathfrak{h}}

\def\cF{\mathcal{F}}

\def\bo{\textbf{o}}

\def\=>{\Longrightarrow}

\def\to{\longrightarrow}

\def\o+{\oplus}
\def\bo+{\bigoplus}

\def\<{\langle}
\def\>{\rangle}
\def\({\left(}
\def\){\right)}

\def\^{\wedge}
\def\+{\dagger}

\def\dd[#1,#2]{\frac{d#1}{d#2}}
\def\del[#1,#2]{\frac{\partial #1}{\partial #2}}
\def\over[#1]{\overline{#1}}
\def\vec[#1]{\overrightarrow{#1}}

\makeatletter
\def\mr@ignsp#1 {\ifx\:#1\@empty\else #1\expandafter\mr@ignsp\fi}%
\newcommand{\multiref}[1]{\begingroup
\xdef\mr@no@sparg{\expandafter\mr@ignsp#1 \: }%
\def\mr@comma{}%
\@for\mr@refs:=\mr@no@sparg\do{\mr@comma\def\mr@comma{,}\ref{\mr@refs}}%
\endgroup}
\makeatother

\newcommand{\hypref}[2]{\ifx\href\asklFrenkel:2013udaas #2\else\href{#1}{#2}\fi}

\tikzset{->-/.style={decoration={
  markings,
  mark=at position .5 with {\arrow{latex}}},postaction={decorate}}}
\tikzset{
    >=latex
    }

\newcommand{\nc}{\newcommand}
\nc{\on}{\operatorname}
\nc{\la}{\lambda}
\nc{\wh}{\widehat}
\nc{\ghat}{\wh\g}
\nc{\mb}{\mathbf}

\begin{document}
\title[Superopers revisited]{Superopers revisited}

\author[A.M. Zeitlin]{Anton M. Zeitlin}
\address{
          Department of Mathematics, 
          Louisiana State University, 
          Baton Rouge, LA 70803, USA\newline
Email: \href{mailto:zeitlin@lsu.edu}{zeitlin@lsu.edu},\newline
 \href{http://math.lsu.edu/~zeitlin}{http://math.lsu.edu/$\sim$zeitlin}
}

\date{\today}

\numberwithin{equation}{section}

\begin{abstract}
The relation between special connections on the projective line, called Miura opers, and the spectra of integrable models of Gaudin type provides an important example of the geometric Langlands correspondence. The possible generalization of that correspondence to simple Lie superalgebras is much less studied. Recently some progress has been made in understanding the spectra of Gaudin models and the corresponding Bethe ansatz equations for some simple Lie superalgebras. At the same time, the original example was reformulated in terms of an intermediate object: Miura-Pl\"ucker oper. It has a direct relation to the so-called $qq$-systems, the functional form of Bethe ansatz, which, in particular, allows $q$-deformation. In this note, we discuss the notion of superoper and relate it to the examples of $qq$-systems for Lie superalgebras, which were recently studied in the context of Bethe ansatz equations. We also briefly discuss the $q$-deformation of these constructions. 
\end{abstract}

\maketitle

\setcounter{tocdepth}{1}
\tableofcontents

\section{Introduction}

One of the most well-understood examples of geometric Langlands correspondence, studied by  E. Frenkel and collaborators \cite{Feigin:1994in}, \cite{Frenkel:1995zp}, \cite{Frenkel:2003qx},\cite{Frenkel:2004qy}, \cite{Feigin:2006xs}, \cite{Rybnikov:2010}, is  the relation between the spectrum of Gaudin models for simple Lie algebra $\mathfrak{g}$ and Miura $^LG$-opers, namely certain meromorphic connections for the principal $^LG$-bundles over $\mathbb{P}^1$. 
Here $^LG$ is a simple Lie group of adjoint type with a Lalnglands dual Lie algebra $^L\mathfrak{g}$. In particular, the set of algebraic equations, known as Bethe equations, which describe the spectrum of the Gaudin model, provides the relation between the ``moduli" parameters determining the corresponding oper connections.

Recently a $q$-deformation of this correspondence was studied \cite{KSZ}, \cite{Frenkel:2020}, which led to various generalizations \cite{KZminors}, \cite{Koroteev:2020tv} and applications \cite{Koroteev:2021a}, \cite{Koroteev2022zoo}.  
In particular, the concepts of $Z$-twisted Miura oper and $Z$-twisted Miura-Pl\"ucker oper were introduced, which have analogues in the differential case \cite{Brinson:2021ww}, \cite{zeitlin2022wronskians}. 
The first notion means that a meromorphic gauge transformation produces a constant connection $Z$ out of a given Miura oper.  
The second, $Z$-twisted Miura-Pl\"ucker oper, is a less restrictive notion, which one can view as a first approximation to a $Z$-twisted condition. Namely, one applies $Z$-twisted condition only to the induced $GL(2)$-oper connections in the 2-dimensional subbundle of the associated bundle for each fundamental representation, corresponding to the vectors of fundamental weight and its elementary Weyl reflection. As a result, one obtains one-to-one correspondence between the data of $Z$-twisted Miura-Pl\" ucker opers and the functional system of equations called $qq$-system ($QQ$-system in $q$-deformed case). Such functional relations previously emerged in the study of the Gaudin model and other spin chain models. In that case, parameter $Z$ is a twist parameter for the boundary condition. Under certain nondegeneracy conditions, these functional relations lead to the Bethe equations describing the spectrum of the related integrable model. 
It was also shown in \cite{Brinson:2021ww} that upon these nondegeneracy conditions, $Z$-twisted Miura-Pl\"ucker opers turn into $Z$-twisted Miura opers, where $Z$ is an element of Cartan. 
An important notion in proving that is the notion of B\"acklund transformations, the gauge transformations, which produce $w(Z)$-twisted Miura oper from a given one, where $w$ is an element of the Weyl group of $^LG$.  
Miura oper connections are specializations of more general oper connections: they are required to preserve specific reduction to the Borel subgroup of the principal bundle to the Borel subgroup. In other words, the B\"acklund transformations allow to ``travel"  between various Miura oper connections corresponding to one single oper connection.  In the context of Bethe ansatz on the level of $qq$-systems, these transformations are known as {\it reproduction procedure} and  {\it populations}, and were tied to the transformations of scalar differential operators, see, e.g., \cite{mukhvar2008}, \cite{MV2002}. 

One would expect to repeat that story in the context of simple superalgebras and supergroups: that, in particular, would give an example of Langlands duality in the super-context. The quest for finding a notion of oper in the context of simple supergroups started in \cite{Zeitlin:2013iya}, where the notion of oper on a super Riemann surface, called {\it superoper}, was introduced, which uses the concept of flat superholomorphic connection. We considered the case of supergroups that allow a purely fermionic system of simple roots. In that case, there is an explicit representation of the space of opers in terms of sections of super projective connections and superconformal vector bundles, known as superconformal fields. Also, as a particular example, we have shown in \cite{Zeitlin:2013iya} the correspondence between the $OSP(1|2)$-superoper on a super Riemann sphere and $\mathfrak{osp}(1|2)$-Gaudin model studied in \cite{kulish2001bethe}, \cite{kulish2001creation}. 

Here we consider the superopers for general simple simply-connected Lie supergroups on a super Riemann sphere. We associate with them a purely even connection on a projective line, which is an oper connection for the reductive purely even subgroup.  We define the notion of $Z$-twisted Miura-Pl\"ucker oper in this case and formulate more general notions of 
$qp$-system and $qq$-systems describing those and their specializations. 
The difference between purely even cases and supergroup cases is that there are several Dynkin diagrams for simple superalgebras, each generally producing a gauge inequivalent oper connection. At the same time, one can ``travel" between various Dynkin diagrams for a given Lie supergroup if one adds odd Weyl reflections generated by odd roots, which are not automorphisms of Lie superalgebra and certainly do not lift to the Lie supergroup.

In this note, we propose to unify the Miura oper connections corresponding to all Dynkin diagrams by introducing formal B\"acklund transformations, which are no longer produced by the gauge transformations of a given connection, which follow the generalized Weyl transformation property on the level of $qq$-systems. We conjecture that satisfying these B\"acklund transformations under certain nondegeneracy conditions will put enough constraints on the $qq$-system to be in 1-to-1 correspondence with Bethe equations for the corresponding Gaudin model. 

Recently, in \cite{huang2019diff}, the authors looked at the $qq$-system, which describes the Bethe ansatz equations for $\mathfrak{gl}(n|m)$ Gaudin model, and studied their populations along the lines of \cite{MV2002}.  That work was continued in paper \cite{lu2021bethe}, where the authors aimed to do the same for $\mathfrak{osp}(n|2m)$ Gaudin Bethe ansatz using the same technique of the transformations of the factorizations of pseudo-differential operators. 
The resulting transformations between data of the populations of these $qq$-systems bring remarkable transformations for the data which characterize the Miura-Pl\"ucker oper.  We hope this more geometric approach can illuminate these reproduction procedures. 

The structure of the paper is as follows. In Section 2 we discuss, following \cite{Zeitlin:2013iya}, the notion of Miura superoper, which a superholomorphic connection \cite{RT} on a super Riemann surface \cite{baranov1985multiloop}, \cite{friedan1986notes}. We also define its reduction to the oper connection on a Riemann surface. In Section 3, we introduce the updated notion of Miura-Pl\"ucker opers and the $qp$-system, which generalizes $qq$-system for opers with regular singularities in pure even case. Section 4 is mainly devoted to the discussion of B\"acklund transformations for Miura opers related to the extended Weyl group and relation to the results of \cite{huang2019diff}, \cite{lu2021bethe}. In the end, we also speculate about the $q$-deformation of these constructions, which should be related to the recent work \cite{huang2019rat}. 

\subsection*{Acknowledgements}
A.M.Z. is partially supported by Simons
Collaboration Grant 578501 and NSF grant
DMS-2203823.

\section{Superopers on supercurves and their reduction}

\subsection{A reminder on super Riemann surfaces and super holomorphic 
connections}
For the general information on supermanifolds and superschemes one can consult \cite{manin2013gauge}, \cite{berezin2013introduction}, \cite{leitessup}. For supercurves and Riemann surfaces one could follow \cite{bergvelt1999supercurves}, \cite{witten2012notes}; for connections for vector bundles over super-Riemann surfaces we refer to \cite{RT}.  

A $supercurve$ of dimension $(1|1)$ over some fixed Grassmann algebra $S$ (which is fixed throughout this paper) is a pair $(X,\mathcal{O}_{X})$, where $X$ is a topological space and $\mathcal{O}_X$ is a sheaf of supercommutative $S$-algebras over $X$ such that $(X,\mathcal{O}^{\rm{red}}_{X})$ is an algebraic curve: $\mathcal{O}^{\rm{red}}_{X}$ is obtained from $\mathcal{O}_X$ by getting rid of nilpotents. 
Locally, for some open sets $U_{\alpha}\subset X$ and some linearly independent elements $\{\theta_{i}\}$ we have $\mathcal{O}_{U_{i}}=\mathcal{O}^{\rm{red}}_{U_{i}}\otimes S[\theta_{i}]$.

 Such collection of open sets $\{U_{i}\}$ serve as coordinate neighborhoods for supercurves with coordinates $(z_{i}, \theta_{i})$.  
The coordinate transformations on the overlaps $U_{i}\cup U_{j}$ 
are given by the following formulas: $z_{i}=F_{ij}(z_{j}, \theta_{j})$, $\theta_{i}=\Phi_{ij}(z_{j}, \theta_{j})$, where $\{F_{ij}\}$, $\{\Phi_{ij}\}$ are even and odd functions correspondingly.

A super Riemann surface $\si$ over some Grassmann algebra $S$ is a supercurve of dimension $1|1$ over $S$, with one more extra structure: there is 
a subbundle $\mathcal{D}$ of $T\Sigma$ of dimension $0|1$, such that for any nonzero section $D$ of $\mathcal{D}$ on an 
open subset $U$ of $\si$, $D^2$ is nowhere proportional to $D$, i.e. one obtains the exact sequence:
\begin{eqnarray}\label{exact}
0\to \mathcal{D}\to T\si\to \mathcal{D}^2\to 0.
\end{eqnarray}
One can pick the holomorphic local coordinates in such a way that this odd vector field 
will have the form $f(z,\theta)D_{\theta}$ for non-vanishing function $f(z,\theta)$ and
\begin{eqnarray}
D_{\theta}=\partial_{\theta}+\theta\partial_z, \quad D_{\theta}^2=\partial_z.
\end{eqnarray}
Such coordinates are called $superconformal$. The transformation between two superconformal coordinate systems 
$(z, \theta)$, $(z', \theta')$ is determined by the condition that $\mathcal{D}$ should be preserved:
\begin{eqnarray}
D_{\theta}=(D_{\theta}\theta') D_{\theta'},
\end{eqnarray}
so that the constraint on the transformation coming from the local change of coordinates is  
$D_{\theta} z'-\theta'D_{\theta} \theta'=0$.
An important  example of a super Riemann surface is the super Riemann sphere $SC^*$: there are two  
charts $(z, \theta)$, $(z, \theta')$ so that
$$
z'=-\frac{1}{z},\quad \theta'=\frac{\theta}{z}.
$$
We call the sections of 
$\mathcal{D}^n$  the {\it superconformal fields} of dimension $-n/2$ following  \cite{witten2012notes}. 
In particular, taking the dual of the exact sequence (\ref{exact}), 
we find that a bundle of superconformal fields of dimension 1, namely $\mathcal{D}^{-2}$,  is a subbundle in $T^*\si$. Considering the superconformal coordinate system, a nonzero section of 
this bundle is generated by $\eta=dz-\theta d\theta$, which is orthogonal to $D_{\theta}$ under the standard pairing.

Let us consider a principal bundle $\mathcal{F}_G$ over the super Riemann surface with the Lie supergroup $G$ over Grassmann algebra $S$. As 
usual, locally one can associate to the connection a differential operator, so that in the chart $(z,\theta)$ the connection has the following form:
\begin{eqnarray}
&&d_A=d+A=d+(\eta A_z+d\theta A_{\theta})+(\bar{\eta}A_{\bar z}+d\bar\theta A_{\bar\theta})=\nonumber\\
&&(\partial+\eta A_z+d\theta A_{\theta})+(\bar \partial+\bar \eta A_{\bar z}+d\bar \theta A_{\bar \theta})=(\eta D^A_z+d\theta D_\theta^A)+(\bar \eta D^A_{\bar z}+d\bar \theta D_{\bar \theta}^A).
\end{eqnarray}  
Here $A$ takes values in $\mathfrak{g}_S$, the Lie algebra of $G$. We note, that we used here the fact that $d=\partial+\bar \partial$ and $\partial=\eta\partial_z+d\theta D_{\theta}$. 
The expression for 
the curvature is:
\begin{eqnarray}
F=d_A^2=d\theta d\theta F_{\theta\theta}+\eta d\theta F_{z\theta}+ d\bar \theta d\bar \theta F_{\bar \theta\bar \theta}+ \bar\eta d\bar \theta 
F_{\bar z\bar \theta}+\eta\bar \eta F_{z\bar z}+ \eta d\bar \theta F_{z\bar \theta}+\bar \eta d\theta F_{\bar z\theta}+d\theta d\bar \theta F_{\theta\bar \theta},\nonumber
\end{eqnarray} 
where 
$$F_{\theta\theta}=- {D^A_{\theta}}^2+D^A_z, ~F_{z\theta}=[D^A_z, D^A_{\theta}], ~F_{z, \bar z}=[D^A_z, D^A_{\bar z}], ~F_{z\bar \theta}=
[D^A_z, D^A_{\bar \theta}], ~F_{\theta\bar \theta}=-[D^A_{\theta}, D^A_{\bar \theta}].$$

It appears that if the connection $d_A$ offers partial flatness, which implies   $F_{\theta\theta}=F_{z\theta}=F_{\bar\theta\bar \theta}=F_{\bar z
\bar \theta}=0$, then there is a superholomorphic structure on any associated bundle (i.e. transition functions of the bundle can be 
made superholomorphic) \cite{RT}. We are interested in the flat superholomorphic connections. In this case, since 
$F_{\theta\theta}=0$, the connection is fully determined by the $D^A_{\theta}$ locally. In other words it is determined by the 
following odd differential operator, which from now on will denote $\hat{\nabla}$: 
\begin{eqnarray}\label{fc}
\hat{\nabla}=D_{\theta}+A_{\theta}(z, \theta), 
\end{eqnarray}
so that the gauge transformation properties  for $A_{\theta}$ are: $A_{\theta}\to gA_{\theta}g^{-1}-D_{\theta}g g^{-1}$, where $g$ is a superholomorphic function providing change of trivialization. \\

\subsection{Miura Superopers}


\subsubsection{Notations}
We refer to \cite{kac1977lie}, \cite{kac2006representations},  \cite{berezin2013introduction}, \cite{cheng2012dualities}, \cite{frappat2000dictionary} for further information regarding simple Lie supergroups, superalgebras, and their representations. Let $G$ is a simple simply connected Lie supergroup of rank $r$
 over some Grassmann algebra $S$, $B_-$ is its fixed Borel subgroup associated to a given Dynkin diagram with unipotent radical $N_-=[B_-, B_-]$. 
 Let $B_+$ be the opposite Borel subgroup containing $H$ and
$N_+=[B_+,B_+]$. 
Note that the Lie algebra $\mathfrak{g}_S$ of ${G}$ is a module over $S$, namely $\mathfrak{g}_S=S\otimes \mathfrak{g}$, where $\mathfrak{g}$ is a simple Lie superalgebra over $\mathbb{C}$. 
Let $\{ \alpha_1,\dots ,\alpha_r \}$ be the set of
positive simple roots for the pair $H\subset B_+$. For a given Dynkin diagram, 
we divide the index set of simple roots 
$I=\{1,\dots , r\}$ into union 
$I=I_w  \sqcup I_g \sqcup I_b$ corresponding to the index set of white, grey, and black roots. $W$ stands for the Weyl group generated the Weyl reflections corresponding to $\{\alpha_i\}_{i\in I_w}$. 

Let $\{e_i, f_i, \check{\alpha}_i\}_{i=1, \dots, r}$ be the Chevalley generators of $\mathfrak{g}$, $a_{ji}=\langle \check{\alpha}_j, \alpha_i\rangle$ is the Cartan matrix. We note also that for grey roots $a_{ii}=0$ for $i\in I_g$, and, 
in additional to standard Serre relations for $\{e_i\}$, $\{f_i\}$ there are extra Serre relations related to the grey root generators.

The Lie superalgebra  $\fb_{-,S}=\fb_{-}\otimes S=\Lie(B_-)$ is generated by the $f_i$'s and the
$\check{\alpha}_i$'s, and  $\fb_{+, S}=\fb_{+}\otimes S=\Lie(B_+)$ is generated by the $e_i$'s and the $\check{\alpha}_i$'s.

Let $\pi: S\to \mathbb{C}$ is the natural projection.  That can be extended to $\pi: G\rightarrow \bar{G}$, where $\bar{G}$ is the underlying reductive simply connected group over $\mathbb{C}$. We denote by 
$\bar{\mathfrak{g}}$, 
$\bar{\mathfrak{b}}_{\pm}$, 
$\bar{\mathfrak{h}}$, 
$\bar{\mathfrak{n}}_{\pm}$ the corresponding pure even reductive Lie algebra, and the pure even versions of its Borel, Cartan and the maximal nilpotent subalgebras, while $\bar{G}$, 
$\bar{B}_{\pm}$, 
$\bar{H}$, 
$\bar{N}_{\pm}$ the corresponding subgroups.

\subsubsection{Definition of superopers}

Now we are ready to define the notion of superoper following a similar definition in the pure even case \cite{Beilinson:2005}, \cite{Frenkel:2003qx} and inspired by the study of integrable hierarchies of Drinfeld-Sokolov type and related integrable models \cite{olshanetsky1983supersymmetric}, \cite{inami1991lie}, \cite{gualzetti1993quantum}, \cite{delduc1998supersymmetric}.

Let us consider a principal $G$-bundle $\mathcal{F}_G$ over a super Riemann surface 
$X$ and its reduction $\mathcal{F}_{B_-}$ to the Borel subgroup $B_-$. 
We assume that it has a flat superholomorphic connection determined by $\hat{\nabla}$.
suppose $\hat{\nabla}'$ is another superholomorphic connection, which preserves $\mathcal{F}_{B_-}$. Then the 
difference $\hat{\nabla}'-\hat{\nabla}$ has a structure of superconformal field of dimension $1/2$ with values in the associated bundle 
$\mathfrak{g}_{\mathcal{F}_{B_-}}$.

Following the purely even case we define an open $B_-$-orbit 
${\bf O}_S\subset[\mathfrak{n}_{-_S}, \mathfrak{n}_{-,S}]^{\perp}/\mathfrak{b}_{-,S}$, consisting of vectors, stabilized by $N_-$ and such that all the simple root components of these vectors with respect to the adjoint action of $H$ are non-zero. Here the orthogonal component is taken with respect to the nondenerate form for a given simple Lie superalgebra. 

\begin{Def}\cite{Zeitlin:2013iya}
 A {\it G-superoper} on a super Riemann surface $\si$ is the triple $(\mathcal{F}, \mathcal{F}_{B_-}, \mathcal{\nabla})$, where $\mathcal{F}$ is a principle 
$G$-bundle, $\mathcal{F}_{B_-}$ is its $B_-$-reduction and $\nabla$ is a long superderivative on $\mathcal{F}$, such that 
$\nabla/\mathcal{F}_{B_-}$ takes values in ${\bf O}_{\mathcal{F}_B}$. 
\end{Def}

Therefore, locally on the open subset $U$, with coordinates $(z, \theta)$, with respect to the 
trivialization of $\mathcal{F}_B$, the structure of the 
superholomorphic connection is:
\begin{eqnarray}\label{sops}
\hat{\nabla}=D_{\theta}+\sum^r_{i=1}a_i(z, \theta)e_i+b(z,\theta),
\end{eqnarray}
where each $a_i(z, \theta)$ is a nonzero function of opposite parity to $e_i$ and $b(z, \theta)$ is an odd $\mathfrak{b}_{-,S}$-valued function. 

\subsubsection{Miura superopers and their pure even counterparts}

Let us start from a definition of  Miura superoper and  

\begin{Def}    \label{Miura}
  A {\em Miura $G$-superoper} on $X$ is a quadruple
 $(\cF_G,\nabla,\cF_{B_-},\cF_{B_+})$, where $(\cF_G,\hat{\nabla},\cF_{B_-})$ is a
  meromorphic $G$-oper on $SC^*$ and $\cF_{B_+}$ is a reduction of
  the $G$-bundle $\cF_G$ to $B_+$ that is preserved by the
  connection $\hat{\nabla}$.
\end{Def}

From now on we set $\si=SC^*$. Then one can put Miura superoper in the following canonical form.

\begin{Prop}    \label{miurastform}
For any Miura $G$-oper on $SC^*$, there exists a
trivialization of the underlying $G$-bundle $\cF_G$ on an open
dense subset of $SC^*$ for which the superoper connection has the form
\begin{equation}    \label{genmiura}
\hat{\nabla}=D_{\theta}+\sum^r_{i=1}g_i(z,\theta)\check{\alpha}_i+\sum^r_{i=1}{a_i(z,\theta)}e_i,
\end{equation}
where $g_i(z,\theta), a_i(z)\in {S}(z)[\theta]$, so that $g_i(z,\theta)$ are all even and $a_i(z,\theta)$ are opposite in parity to $e_i$.
\end{Prop}

The proof of this proposition is similar to the one in \cite{Brinson:2021ww} with the use of the cell partition via super extension of Weyl group in \cite{manvor2, manvor1}, see also \cite{penkov1990borel}.

From now on we assume that $\{D_{\theta}a_i(z,\theta)\}_{i\in I_w}$, as well as  $\{a_i(z,\theta)_{i\in I_{g}\cup I_{b}}$ are invertible.

\begin{Def}    \label{Ztwoper}
  A {\em ${\hat{Z}}$-twisted $G$-superoper} on $SC^*$ is a $G$-superoper
  that is equivalent to the constant element $A(\theta, z)=\theta {\hat Z}$, where $\hat Z \in \mathfrak{g}_S \subset \mathfrak{g}(z,\theta)$ under the gauge action of $G(z, \theta)$.
\end{Def}

For simplicity from now on we will assume that $\hat Z\in \mathfrak{h}_S$ is regular semisimple. One can generalize most of the results to $Z\in\mathfrak{b}_{+,S}$ as it was done in \cite{Brinson:2021ww}.

\begin{Rem}
Note, that instead one could consider any element $Z'=\xi+\theta Z$, where $xi$ i an odd  element of $\mathfrak{g}_S$, instead of $\theta Z$ in the above definition, but one can remove $xi$ by gauge transformation given by $\exp(\theta\xi)\in G(\theta)$.  
\end{Rem}

Now we want to get rid of extra odd variables to see the relation of superopers to opers on $\mathbb{P}^1$ for a certain reductive group.

First, we will get rid of the $\theta$ variable.  Let us represent the superoper connection as 
\begin{eqnarray}
\hat{\nabla}=D_{\theta}+\theta M(z) +N(z)
\end{eqnarray}
making the dependence on $\theta$ explicit. The gauge transformations can be factorized this way: $$g(z,\theta)=(1-\theta R(z))U(z),$$ 
where $R(z)\in \mathfrak{g}(z)$, $U(z)\in G(z)$. There is a unique $R(z)$ , namely  $R(z)=N(z)$, such that 
\begin{eqnarray}
g(z,\theta)^{-1}\hat{\nabla}g(z, \theta)=\partial_{\theta}+\theta U^{-1}(z) \tilde{\nabla}U(z),
\end{eqnarray}
where the meromorphic connection $\tilde{\nabla}$ is given locally by a differential operator
\begin{eqnarray}
\tilde{\nabla}=\partial_z+\frac{1}{2}[N(z),N(z)]+M(z). 
\end{eqnarray}

Thus we obtain the following Proposition.

\begin{Cor} $Z$-twisted condition for Miura $G$-superoper just implies that connection 
$\tilde{\nabla}$ on $\mathbb{P}^1$ is gauge equivalent to a constant connection $\hat{Z}\in \mathfrak{g}$.
\end{Cor}

Now let us proceed to Miura $G$-superoper, so that the resulting connection is as in (\ref{genmiura}) and let's construct such connection 
$\tilde\nabla$. This way, we obtain $G$-connection on $\mathbb{P}^1$. 
Now, let us apply the map $\pi: S\rightarrow \mathbb{C}$,  which strips dependence on all the odd parameters. This way, we obtain the connection $\nabla=\pi (\tilde{\nabla})$ on a principal $\bar{G}$-bundle, which have locally the following form:
\begin{eqnarray}\label{paramop}
\nabla\equiv \pi(\tilde{\nabla})=\partial_z+u(z)+\sum_{i\in I_w} L_i(z)e_i+\frac{1}{2}\sum_{i,j\in I_b\cup I_g } L_i(z)L_j(z)[e_i, e_j],
\end{eqnarray} 
where $\{L_i(z)\}_{i=1, \dots, r}$ are nonzero  rational functions, 
so that in the original connection (\ref{genmiura}):
\begin{eqnarray}
\pi (a_i(z,\theta))=
\left\{ 
  \begin{array}{ c l }
   L_i(z) & \quad \textrm{if} ~ i\in I_w \\
      \theta L_i(z)             & \quad \textrm{if} ~i\in I_b\cup I_g.
  \end{array}
\right.
\end{eqnarray}


That gives rise to the following definition.

\begin{Def}
We say that the quadruple 
$(\mathcal{F}_{\bar{G}}, \mathcal{F}_{{\bar B}_+},  \mathcal{F}_{\bar B_-},\nabla)$, 
where $\nabla$ is a connection on a principal bundle
 $\mathcal{F}_{\bar G}$ over 
 $\mathbb{P}^1$ together with reductions $\mathcal{F}_{{\bar B}_\pm}$ to Borel subgroups ${\bar B}_{\pm}$ is a Miura oper associated to Miura G-superoper  
$(\mathcal{F}_G, \mathcal{F}_{B_+},  \mathcal{F}_{B_-},\hat{\nabla})$.
\end{Def}

In particualr, we notice that the differential operator which isolates the Cartan part of \ref{paramop}: 
\begin{eqnarray}
\nabla^H=\partial_z+u(z)
\end{eqnarray}
defines an ${H}$-connection on $\mathbb{P}^1$, which we call a {\it Cartan connection } 
$\nabla^H$ associated to Miura oper $\nabla$. We say that $H$-connection is $Z$-twisted if it is gauge equivalent to the constant connection $\partial_z+Z$, where $Z\in \mathfrak{h}$.

\begin{Rem}
i)We note here that in general, even for distinguished Dynkin diagrams the collection  
 $\{\check{\alpha_i}\}_{i=1, \dots, r}$, $[e_i,e_j]$ for all $i,j\in I_g\cup I_b$ and $e_i$ for all $i\in I_w$ do not give Chevalley generators producing $\bar{B}_+$, but instead a  Borel subgroup of a smaller reductive group.\\
ii) If the superoper was $Z$-twisted, he resulting connection $\nabla$ is gauge equivalent to $\pi(Z)$.
\end{Rem}

\section{Miura-Pl\"ucker opers,  $qp$-systems, and $qq$-systems}

\subsection{$Z$-twisted Cartan connections} Let us fix the component notation for the Cartan connection, associated to a given superoper and the corresponding twist element $Z$:
\begin{eqnarray}
u(z)=\sum_{i\in I_w+I_g}u^i(z)\check{\alpha}_i+\sum_{i\in I_b}u^i(z)\frac{\check{\alpha}_i}{2}; \quad Z=\sum_{i\in I_w\cup I_g}\zeta_i\check{\alpha}_i+\sum_{i\in I_b}\zeta_i\frac{\check{\alpha}_i}{2} 
\end{eqnarray}

This leads to the following proposition.

\begin{Prop}\label{cartztwist}
If the Cartan connection parametrized by $u(z)$ in (\ref{paramop}) is $Z$-twisted, then there exist rational functions $\{p_i(z)\}_{i=1, \dots, r}$ such that 
\begin{equation}
u_i(z)=\zeta_i+ \ln'\big[p_i(z)\big], \quad 1, \dots, r.
\end{equation}
\end{Prop}

One can view it as a first approximation for $Z$-twisted condition for the Miura $\bar{G}$-oper connection (\ref{paramop}).

Now let us consider in detail the rank $r=1$ examples of $Z$-twisted Miura opers.

\subsection{Low rank cases}
\subsubsection{Z-twisted Miura $SL(2)$-opers}

In this case we are dealing with the meromorphic $SL(2)$-connection for $SL(2)$ bundle 
which has the following form:
\begin{equation}\label{miurasl2}
\nabla=\partial_z+u(z)\check{\alpha}+L(z)e
\end{equation}
where $e, f, \check{\alpha}$ are the Chevalley generators of $\mathfrak{sl}(2)$. The $Z$-twisted condition states that there exists $U(z)\in B_+(z)$, such that 
\begin{equation}\label{sl2ztw}
U(z)\nabla U^{-1}(z)=\partial_z+Z. 
\end{equation}
where $Z=\zeta \check{\alpha}$. 
We can represent the resulting group element as 
$$
U(z)=e^{q(z)e} {p(z)}^{\check{\alpha}},
$$
where $p(z), q(z)\in \mathbb{C}(z)$. Looking at the coefficient of $\check{\alpha}$ in (\ref{sl2ztw}), we obtain 
\begin{equation}\label{ueq}
u(z)=\zeta+\ln '\big[p(z)\big] 
\end{equation}
The $e$-coefficient gives the equation:
$$
q'(z)+2\zeta q(z)=p^2(z)L(z)
$$
Representing 
$q(z)=\frac{q_-(z)}{q_{+}(z)}$, 
so that 
$q_{\pm}(z)\in \mathbb{C}[z]$, $q_+(z)$ is monic, we obtain the following equations:
\begin{equation}\label{wlambda}
 W(q_-,q_+)(z)+2\zeta q_-(z)q_+(z)=\Lambda(z), 
\end{equation}
so that 
\begin{equation}
\Lambda(z)=q_{+}^2(z)p^2(z)L(z)
\end{equation}
is a polynomial.
Expanding $u(z)=\frac{\tilde{u}(z)}{\tilde{\Lambda}(z)}$ and  
$p(z)=\frac{p_-(z)}{p_+(z)}$, where $p_{\pm}(z)$ do not have common roots and $p_{-}(z)$ is chosen to be monic, we obtain:
\begin{equation}\label{wlambdat}
\tilde{\Lambda}(z)=p_{+}(z)p_{-}(z).
\end{equation}
Note, that given the factorization of denominator $\tilde{\Lambda}$, the numerator 
is determined uniquely:
$$\tilde{u}(z)=W(p_-,p_+)(z)+2\zeta p_-(z)p_+(z).$$

Thus let us call the following system of equations:
\begin{eqnarray}\label{sl2pq}
&& W(q_-,q_+)(z)+2\zeta q_-(z)q_+(z)=\frac{q_{+}^2(z)p_{-}^2(z)L(z)}{p_+(z)^2},\\
&& p_{+}(z)p_{-}(z)=\tilde{\Lambda}(z)\nonumber
\end{eqnarray}
the {\it $pq$-system} for $\mathfrak{sl}(2)$.
\begin{Prop}
There is one-to-one correspondence between $Z$-twisted Miura $SL(2)$-opers with the connection (\ref{miurasl2}) and solutions of the $pq$-system (\ref{sl2pq}), where $p(z)=\frac{p_{-}(z)}{p_{+}(z)}$  as well as $q(z)=\frac{q_{-}(z)}{q_{+}(z)}$ are irreducible fractions, so that  $\tilde{u}(z)=W(p_-,p_+)(z)+2\zeta p_-(z)p_+(z)$.
\end{Prop}
A simplification of the $\mathfrak{sl}(2)$ $pq$-system is given by the following identification
\begin{eqnarray}
\tilde{\Lambda}(z)=p_{+}(z)=q_{+}(z),
\end{eqnarray}
which leaves just one equation:
\begin{eqnarray}
W(q_-,q_+)(z)+2\zeta q_-(z)q_+(z)=\Lambda(z),
\end{eqnarray} 
so that 
$u(z)=\zeta-\ln'(q_+(z))$ and $L(z)=\Lambda(z)$ is a polynomial. That is known as $\mathfrak{sl}(2)$ $qq$-system, which is in one-to one correspondence with the $SL(2)$-opers with regular singularities: the positions of the singularities on $\mathbb{P}^1$ are given by the roots polynomial $\Lambda(z)$. 
Under nondegeneracy conditions, namely $q_{+}(z)$ have distinct roots and have no common roots with $\Lambda(z)$, a simple calculation shows that there is a 
bijection between such nondegenerate solutions of $\mathfrak{sl}(2)$ $qq$-system and the the Bethe equations of $\mathfrak{sl}(2)$ Gaudin model:
\begin{eqnarray}
&&2\zeta+\partial_z\log\Big[\Lambda(z)(z-w_\ell)^2\Big]\Bigg|_{z=w^i_\ell}=0,\\ 
&& i=1,\dots, r; \quad  \ell=1, \dots, \deg(q_+(z)).\nonumber
\end{eqnarray}

\subsubsection{Z-twisted Miura $SL(1|1)$ opers and abelian connections}
 Although it is not a simple superalgebra because of nontrivial center, one can apply the notion of oper to $SL(1|1)$-group, which technically corresponds to the grey node of Dynkin diagram, and  it is still useful to consider it.   
We see that in this case the connection $\bar{\nabla}$  and $\nabla$ are reduced to the abelian connection corresponding to the central element $\check{\alpha}$ of $SL(1|1)$:
\begin{equation}\label{mopsl11}
\nabla=\partial_z+u(z)\check{\alpha}.
\end{equation}
The $Z$-twisted condition leads to the equation (\ref{ueq}). Expressing $u(z)=\frac{\tilde{u}(z)}{\tilde\Lambda(z)}$, 
one can say that this is a particular case of the first example when $L(z)=0$, so one can  call the equation 
\begin{eqnarray}\label{pqsl11}
p_+(z)p_-(z)=\tilde{\Lambda}(z)
\end{eqnarray}
the $\mathfrak{sl}(1|1)$ $pq$-system, although $q$-part is absent here.
\begin{Prop}
There is one-to-one correspondence between $SL(1|1)$ -opers (\ref{mopsl11})  and the solutions to the equation (\ref{pqsl11}), so that $p_{\pm}(z)$ 
have no common roots and  
$${u}(z)=\tilde{\Lambda}(z)^{-1}\Big[W(p_-,p_+)(z)+2\zeta p_-(z)p_+(z)\Big].$$

\end{Prop}
The equation $(\ref{pqsl11})$ actually gives the Bethe ansatz solutions for $\mathfrak{gl}(1|1)$ Gaudin model. Namely, let's suppose we can re-express  
\begin{equation}
\tilde{\Lambda}(z)=\ln'(\Lambda(z))\pi(z),
\end{equation}
where $\pi(z)=\prod^n_{k=1}(z-z_k)$ where $z_k$ are all distinct roots of $\Lambda(z)=\prod^n_{k=1}(z-z_k)^{d_k}$. 
In this case the $\Lambda$ determines the weights for the appropriate $\mathfrak{gl}(1|1)$ representation and the equation $\ref{pqsl11}$ is equivalent to the Bethe ansatz equations for $\mathfrak{gl}(1|1)$ Gaudin model \cite{mukhin2015gaudin}, \cite{huang2019diff}:
\begin{equation}
\sum^n_{k=1}\frac{d_k}{w_j-z_k}=0,\quad   j=1, \dots \deg(p_+(z)),
\end{equation}
where $\{w_j\}$ are the roots of $p_+(z)$.

\subsubsection{Z-twisted Miura $OSP(1|2)$-opers}
Consider an $\mathfrak {osp}(1|2)$-triple: $e, f, \check{\alpha}$, where $e, f$ are odd Chevalley generators. The corresponding $\overline{OSP(1|2)}$-Miura oper connection is:
\begin{eqnarray}
\nabla=\partial_z + \frac{\check{\alpha}}{2}u(z)+\frac{1}{2}L^2(z)[e,e].
\end{eqnarray}
Notice that this is exactly the Miura $SL(2)$-oper we considered above, since $\frac{\check{\alpha}}{2}$, $[e,e]$ are Chevalley generators of $\mathfrak{sl}(2)$ subalgebra. The only difference is that we have a square as a coefficient of $e^2$. 

Thus we have the following Proposition.

\begin{Prop}
There is a one-to-one correspondence between $Z$-twisted $\overline{OSP(1|2)}$-opers and $Z$-twisted Miura $SL(2)$ opers where $\Lambda(z)=L^2(z)$, $L(z)\in \mathbb{C}[z]$.
\end{Prop}

That correspondence was first noted in \cite{Zeitlin:2013iya}: in particular, it was shown that the resulting Bethe ansatz equations describing coincide with the Bethe ansatz equations for $\mathfrak{osp}(1|2)$ Gaudin model \cite{kulish2001bethe}, \cite{kulish2001creation}.

\subsection{Miura-Pl\"ucker opers}
$Z$-twisted condition is a complicated one to solve. Instead one can look at the  intermediate object. We already have seen the first approximation to that condition given by the Proposition \ref{cartztwist}. Now, we introduce a useful object,  known as $Z$-twisted Miura-Pl\"ucker $\bar{G}$-oper, which is a next interation approximating $Z$-twisted condition. In pure even case, for Miura $G$-opers with regular singularities, it was introduced in \cite{Brinson:2021ww} following the $q$-deformed version in \cite{Frenkel:2020}. 

Let us consider the induced Miura $\bar{G}$-oper $\bar{B}_+$-bundle connections on $\mathcal{V}_i$: the associated bundles, corresponding to highest weight irreducible modules of $\bar{\mathfrak{g}}$: $i)$ $V_{\omega_i}$ if  $i\in I_w\cup I_g$; $ii)$ $V_{2\omega_i}$ when $i\in I_b$. 

Let us define a $B_+$-subbundle $\mathcal{W}_i$, the rank of which will depend on  $i$.
\begin{itemize}
\item If $i\in I_w$, $W_i$ is spanned by the line subbundles 
$\mathcal{L}_i, \tilde{\mathcal{L}}_i$,  which correspond to the vectors of weight $\omega_i$, $\omega_i-\alpha_i$.  
\item If $i\in I_b$, let $\mathcal{W}_i$ be spanned by 
$\mathcal{L}_i$, ${\tilde{\mathcal{L}}}_i$ corresponding to the vectors of weight $2\omega_i$, $2\omega_i-2\alpha_i$. 
\item If $i\in I_g$ we take $W_i$  to be a line bundle  which correspond to the vector of highest weight $\omega_i$. 
\end{itemize}
 Let $\nabla_i$ be the induced connection on $W_i$. 
\begin{Def}\label{defmp}
We say that Miura $\bar{G}$-oper is $Z$-twisted Miura-Pl\"ucker if there exists $v(z)\in B_+(z)$ such that 
\begin{eqnarray}\label{condMP}
\nabla_i=v(z)(\partial_z+Z)v(z)^{-1}|_{W_i}=v_i(z)(\partial_z +Z_i) v_i(z)^{-1},
\end{eqnarray}
where $v_i(z) = v(z)|_{W_i}$ and $Z_i = Z|_{W_i}$.
\end{Def}
The element $v(z)$ is not uniquely determined by Miura-Pl\"ucker oper. Let $\tilde{N}_+(z)$ be the subgroup, generated by all even commutators of $[e_i, e_j]$, $i\neq j$.   
We have the following proposition, which gives equivalence classes of such $v(z)$.
\begin{Prop}
For a given $v(z)$ from the Definition \ref{defmp} any element of coset 
$$v(z)H\tilde{N}_+(z)$$ 
also satisfies (\ref{condMP}). 
\end{Prop}
Following \cite{Brinson:2021ww}, we call such a coset a {\it framing} of Miura-Pl\"ucker oper.

\begin{Def}
The {\it Miura-Pl\"ucker datum} is a pair $(\nabla, v(z)\tilde{N}_+(z))$ consisting of Miura-Pl\"ucker oper and related framing.
\end{Def}

One can fix the corresponding representative in the coset as follows:
\begin{eqnarray}
v(z)=\prod_{i\in I_w\cup I_g}p^i(z)^{-\check{\alpha}_i}
\prod_{j\in I_b}p^j(z)^{-\frac{\check{\alpha}_j}{2}}
\prod_{i\in I_w}
e^{-q^i(z)e_i}
\prod_{j\in I_b}e^{-\frac{1}{2}q^j(z)[e_j,e_j]}. 
\end{eqnarray}

Now we will explore the condition \ref{condMP}, 
For the Cartan part of $\nabla$ we obtain the following equations:
\begin{equation}    \label{uiz}
u^i(z)=\zeta_i+\ln'(p^i(z)).
\end{equation}

Now we will show how it works off-diagonal in each of the cases:

\begin{enumerate}

\item Let $i\in I_w$. Then we have the following condition: 

We first compute the matrix of $v(z)$ and $Z$ acting on the
two-dimensional subspace $W_i$.  The following calculation gives
\begin{equation}    \label{vzz}
v(z)|_{W^i}=
\begin{pmatrix}
  p^i(z)^{-1} & 0\\
  0& p^i(z)\prod_{j\neq i, j\in I_w\cup I_g} p^j(z)^{a_{ji}}\prod_{j\neq i, j\in I_b} p^j(z)^{a_{ji}/2}
 \end{pmatrix}
 \begin{pmatrix}
1 & - \frac{q^i_{-}(z)}{q^i_{+}(z)}\\
 0& 1
 \end{pmatrix}
\end{equation}
and 
\begin{equation}
Z^H|_{W_i}=\begin{pmatrix}
\zeta_i & 0\\
  0& -\zeta_i-\sum_{j\neq i}{a_{ji}}\zeta_j
   \end{pmatrix}.
\end{equation}
That implies the following equation from the top right corner of $2\times 2 $ block:

\begin{equation}    \label{qqw}
\partial_z q^i(z)+\langle Z, \alpha_i\rangle q^i(z)=L_i(z)\Big[{p^i(z)}\Big]^{2}
\prod_{j\neq i, j\in I_w\cup I_g} p^j(z)^{a_{ji}}\prod_{j\neq i, j\in I_b} p^j(z)^{a_{ji}/2}.
\end{equation}

\item If $i\in I_g$, we do not have any extra equations in addition to  (\ref{uiz}).

\item If $i\in I_b$, we are dealing with the same $2\times 2$ block as in $i\in I_w$. Thus we have the following equation:
\begin{equation}    \label{qqb}
\partial_z q^i(z)+\langle Z, \alpha_i\rangle q^i(z)=\Big[L_i(z){p^i(z)}\Big]^{2}
\prod_{j\neq i, j\in I_w\cup I_g} p^j(z)^{2a_{ji}}\prod_{j\neq i, j\in I_b} p^j(z)^{a_{ji}}.
\end{equation}

\end{enumerate}

\subsection{$qp$-and $qq$-system for Lie superalgebra $\mathfrak{g}$.}

\subsubsection {Definition of the $pq$-system and relation to $Z$-twisted Miura-Pl\"ucker opers.}
Let us choose the simple root system of Lie superalgebra $\mathfrak{g}$ with the index set $I=I_w\cup I_b\cup I_g$ and 
the datum of rational functions $\{L_i\}_{i\in I}(z)$, polynomial functions $\{\tilde{\Lambda}_i\}_{i\in I}$ and rational functions $\{p^i(z), q^i(z)\}_{i\in I}$ be the irreducible  fractions:
\begin{equation}
p^i(z)=\frac{p^i_{-}(z)}{p^i_{+}(z)}, \quad q^i(z)=\frac{q^i_{-}(z)}{q^i_{+}(z)}, \quad  i\in I.
\end{equation}
We call the following system of equations: 
\begin{eqnarray}
&&\partial_z q^i(z)+\langle Z, \alpha_i\rangle q^i(z)=F_i(z), \quad i\in I_w\cup I_b,\\\nonumber
&&F_i(z)=\left\{ 
  \begin{array}{ c l }
   L_i(z)\Big[{p^i(z)}\Big]^{2}
\prod_{j\neq i, j\in I_w\cup I_g} p^j(z)^{a_{ji}}\prod_{ j\in I_b} p^j(z)^{a_{ji}/2} & \quad \textrm{if} ~ i\in I_w \\
      \Big[L_i(z){p^i(z)}\Big]^{2}
\prod_{ j\in I_w\cup I_g} p^j(z)^{2a_{ji}}\prod_{j\neq i, j\in I_b} p^j(z)^{a_{ji}}              & \quad \textrm{if} ~i\in I_b
  \end{array}
\right.\\
&&p^i_+(z)p^i_-(z)=\tilde{\Lambda}_i(z), \quad i=1, \dots, r,\nonumber
\end{eqnarray}
the {\it $pq$-system} associated to superalgebra $\mathfrak{g}$.

Using this definition we can restate the result of the previous section as follows.

\begin{Thm}
There is a bijection between the datum of $Z$-twisted Miura-Pl\"ucker opers modulo the data, solutions of the generalized $pq$-systems, so that 
$q^i_{+}(z),q^i_{-}(z)$  as well as $p^i_{+}(z),p^i_{-}(z)$ for all $i=1,\dots, r$ have no common roots and 
\begin{eqnarray}
u_i(z)=\frac{\tilde{u}_i(z)}{\tilde{\Lambda}_i(z)}=\zeta_i+\ln'(p_i)(z).
\end{eqnarray}
\end{Thm}
 
\begin{Rem} 
 We remark here, that we do not make any assumptions/nondegeneracy conditions on far on the roots/poles of the data $\{\tilde{\Lambda}_i\}$, $\{L_i(z)\}$. 
 In the following we will see examples when $\{\tilde{\Lambda}_i(z)\}_{i\in I_g}$ depends on $q^i(z)$. Notice also that there are no equations on $\{L_i(z)\}_{i\in I_g}$.  Those will be take into account when the full $Z$-twisted condition is implemented.
 \end{Rem}
 
 \subsubsection{Reduction to the qq-system, nondegeneracy conditions and Bethe equations}
 
To get in touch with the Bethe ansatz equations for Gaudin model, we will impose the following condition  on the $qp$-system: we require that the reduction of the $pq$-system to the simple even subgroups of $\bar{\mathfrak{g}}$ reproduces the $qq$-system for Miura opers with regular singularities \cite{Brinson:2021ww}. Namely, we assume:
  \begin{equation}
 \tilde{\Lambda}_i(z)=p_+^i(z)=q_+^i(z)~ {\rm  for} ~i\in I_w\cup I_b. 
 \end{equation} 
Also, we redefine: 
\begin{eqnarray}
p^i_{\pm}(z)\equiv q^{i}_{\pm}(z), \quad {\rm for} \quad i\in I_g.
\end{eqnarray}
We call  the resulting system of equations on the $pq$-system data, the $qq$-system associated to $\bar{\mathfrak{g}}$.
\begin{eqnarray}\label{redqq1}
&&W(q^i_-, q^i_+)(z)+\langle Z, \alpha_i\rangle q_+^i(z)q_-^i(z)=F_i(z), \quad i\in I_w\cup I_b,\\\nonumber
F_i(z)&=&\left\{ 
  \begin{array}{ c l }
   L_i(z)
\prod_{j\neq i, j\in I_w} q_+^j(z)^{-a_{ji}}\prod_{ j\in I_b} q_+^j(z)^{-a_{ji}/2} \prod_{ j\in I_g} q^j(z)^{a_{ji}}& \quad \textrm{if} ~ i\in I_w \\
      \Big[L_i(z)\Big]^{2}
\prod_{ j\in I_w} q_+^j(z)^{-2a_{ji}}\prod_{j\neq i, j\in I_b} q_+^j(z)^{-a_{ji}}             
\prod_{ j\in I_g} q^j(z)^{2a_{ji}}  & \quad \textrm{if} ~i\in I_b
  \end{array}
\right.\nonumber\\
&&q^i_+(z)q^i_-(z)=\tilde{\Lambda}_i(z), \quad i\in I_g.\label{redqq2}
\end{eqnarray}
This way we reduced the number of independent functions to $\{L_i(z)\}_{i\in I_w\cup I_b}$ and $\{\tilde{\Lambda}_i(z)\}_{i\in I_g}$. 
Under certain conditions, one can write solutions to the (\ref{redqq1}) in terms of algebraic equations. First, we introduce the nondegeneracy conditions:
\begin{Def}
The solutions to the $qq$-system is called {\it nondegenerate}, if $F_i(z)$ has no common roots with $q^i_+(z)$ and all the roots 
 of $q^i_{+}(z)$ are distinct for all $i$. Moreover $p_+^i(z)$ and $p_-^i(z)$ have no common roots. 
\end{Def}

Then  the following Proposition is true.

\begin{Prop}
If $Z$ is regular semisimple, there is a bijection between nondegenerate solutions to (\ref{redqq1}) and the following algebraic equations
for the roots $\{w^i_{\ell}\}_
{\ell=1, \dots, \deg(q^i_+(z))}$ of $q^i_{+}(z)$:
\begin{equation}\label{bethe1}\begin{gathered}
\langle\alpha_i,Z\rangle+\partial_z\log\Big[F_i(z)(z-w^i_\ell)^2\Big]\Bigg|_{z=w^i_\ell}=0,\\ 
 i=1,\dots, r; \quad  \ell=1, \dots, \deg(q^i_+(z)).
\end{gathered}
\end{equation}

\end{Prop}
The equations (\ref{bethe1}) are known in particular cases as Bethe equations for Gaudin model associated with simple Lie algebras, we will refer to them as Bethe equations of {\it even type}. 
The algebraic equations emerging from (\ref{redqq2}), which are solved just by division  are Bethe equations of {\it odd type}. Of course, in the case, when $\mathfrak{g}$ is a simple algebra, one has only equations of even type. The following proposition follows.   

In the case when $\mathfrak{g}$ is a simple Lie algebra, i.e. $I=I_{w}$ the $qq$-system (\ref{redqq1}),(\ref{redqq2}) reduces to well-known $qq$-system from \cite{mukhvar2008}, \cite{Brinson:2021ww}, imposing the condition that $L_i(z)=\Lambda_i(z)$ is a polynomial:
\begin{eqnarray}
W(q^i_-, q^i_+)(z)+\langle Z, \alpha_i\rangle q_+^i(z)q_-^i(z)=\Lambda_i(z)
\prod_{j\neq i} q_+^j(z)^{-a_{ji}}, \quad i=1, \dots, r,
\end{eqnarray}
where $a_{ji}$ is a Cartan matrix of $\mathfrak{g}$.

The corresponding  $G$-oper connections are called Miura $G$-opers with regular singularities, where the position of singularities are regulated by polynomials $\{\Lambda_i(z)\}_{i=1,\dots, r}$ \cite{Brinson:2021ww}. Locally the connection has the form
\begin{equation}
\nabla=\partial_z+Z-\sum^r_{i=1}\ln'\Big[q^i_+(z)\Big]\check{\alpha}_i +\sum^{r}_{i=1}\Lambda_i(z)e_i.
\end{equation}
\section{Z-twisted Miura opers and the extended Weyl group} 

\subsection{Overview of the pure even case}
In \cite{Brinson:2021ww}, in the case of simple Lie algebras for opers with regular singularities, we have shown that under the nondegeneracy conditions $Z$-twisted Miura-Pl\"ucker opers turn out to be $Z$-twisted Miura opers, i.e., the solutions of the $qq$-system completely determine $Z$-twisted Miura opers. 

Let us look at the related $Z$-twisted oper for regular semisimple $Z$. We find that there are precisely $|W|$ (W is a Weyl group of $\mathfrak{g}$) related $Z$-twisted Miura opers, when $Z$ is  regular semisimple, each described by the solutions of the corresponding $qq$-systems. 
The action of the Weyl group on the space of solutions of such $qq$-systems is given using the following transformations, corresponding to elementary Weyl reflections $w_i$:
\begin{eqnarray}\label{weyl}
Z\to w_i(Z) , \quad  q_{\pm}^j(z)\to 
\left\{ 
  \begin{array}{ c l }
 q^i_{\mp}(z)& \textrm {if} j=i\\
 q^j_{\pm}(z)& \textrm {if} j\neq i
  \end{array}
\right..
\end{eqnarray}
On the level of Miura $G$-oper connections that can be achieved by special gauge transformations from $B_-(z)$:
\begin{eqnarray}\label{eq:PropDef}
\nabla\to e^{\mu_i(z)f_i}~\nabla ~ e^{-\mu_i(z)f_i}, \quad 
\mu_i(z)=\Lambda_i(z)^{-1}\Bigg[\partial_z\log\Bigg(\frac{q^i_-(z)}{q^i_+(z)}\Bigg)+\langle
    \alpha_i,Z\rangle\Bigg],
\end{eqnarray}
 which we called {\it B\"acklund transformations} in \cite{Frenkel:2020}, \cite{Brinson:2021ww}. 
These transformations were previously discussed in the context of Bethe ansatz equations \cite{mukhvar2008}, \cite{mukhvarmiura} leading to the so-called ``populations" of Bethe ansatz equations.
\subsection{Conjectures for simple superalgebras}
In the case of superalgebra $\mathfrak{g}$, Weyl reflections of even roots generate the Weyl group $W$. One can construct a larger group $\tilde{W}$ by adding the reflections, corresponding to the odd roots, which changes the Dynkin diagram for $\mathfrak{g}$  \cite{leites1985embeddings}, \cite{frappat2000dictionary}. In particular, applying such reflections to a given system of simple roots, one can generate all systems of simple roots for a given superalgebra $\mathfrak{g}$. However, these are not automorphisms of $\mathfrak{g}$; of course, one cannot lift these transformations to $G$. 

The relation between $Z$-twisted Miura opers corresponding to a given oper in the pure even case, as discussed in the previous subsection, motivates introducing the following notion, generalizing the one from \cite{Brinson:2021ww}.

\begin{Def}
Consider two $qq$-systems based on Dynkin diagrams based on 
 simple root systems related by a simple reflection $s_i\in \tilde{W}$. 
We say that two solutions of such $qq$-systems are $i$-{\it composable} if they are obtained from each other by 
the transformation (\ref{weyl}) accompanied with certain transformations 
$\{L_i\}, \{\tilde{\Lambda}_i\}\rightarrow \{L^{w_i}_i\}, \{{\tilde \Lambda}^{w_i}_i\}$. We call two solutions of $qq$-systems {\it $w$-composable}, where $w\in \tilde{W}$ if a sequence of such transformations relates them. We call the related $Z$-twisted Miura-Pl\"ucker opers {\it $w$-composable} if their datum expressed via the solution of the $qq$-system is $w$-composable.
\end{Def}

Notice, that we still did not specify the transformations of 
$\{L_i(z)\}_{i\in I_b\cup I_w}$, $\{\tilde{\Lambda}_i(z)\}_{i=1, \dots, r}$, for the $w$-composable $qq$-systems because they may vary from one $qq$-system to another, unlike the pure even case. Also, those may depend on $\{q^j_{\pm}(z)\}_{i=1,\dots, r}$ if $j$ and $i$ are adjacent on the Dynkin diagram. 
In fact, we will see below such an example of $w$-composable family of $qq$-systems associated with $\mathfrak{sl}(n|m)$.
  
Let us formulate the following conjecture, which is the analogue of the main results of \cite{Brinson:2021ww}.

\begin{Con}
i) Under certain nondegeneracy conditions, $Z$-twisted Miura-Pl\"ucker opers are $Z$-twisted Miura opers for certain choices of $\{L_i(z)\}_{i\in I_g}$.  \\
ii) There exists such data of $\{L^w_i(z)\}_{i\in I}$ and $\{\tilde{\Lambda}^w_i(z)\}_{i\in I}$ for the $qq$-systems associated to $\mathfrak{g}$ such that there exists a family of $w$-composable $Z$-twisted Miura opers, which are described by the Bethe equations of Gaudin model for a certain simple superalgebra $^L\mathfrak{g}$.
\end{Con}    
Indeed, this way, we relate various $Z$-twisted opers for different Dynkin diagrams making an object, which one may call a {\it generalized} $Z$-twisted superoper, which unifies all $w$-composable $Z$-twisted Miura opers. 
We do not know what $^L\mathfrak{g}$ could be for the cases beyond $\mathfrak{g}$= $\mathfrak{sl}(m|n)$: in the next subsection we will discuss the discovered class of such $w$-composable $qq$-systems. 

\subsection{What is known: $qq$-systems for $\mathfrak{sl}(n|m)$ and Gaudin models}
In \cite{huang2019rat} a certain version of  $qq$-system was considered in the case of $\mathfrak{sl}(n|m)$ and $Z=0$. For a given Dynkin diagram for $\mathfrak{g}$ one has 
\begin{eqnarray}
&&Wr(q^i_-,q^i_+)(z)=\Lambda_i(z)q_+^{i+1}(z)q^{i-1}_+(z) ~{\rm if}~  i\in I_w,\\\nonumber
&& q^i_-(z)q^i_+(z)=\tilde{\Lambda}_i(z), ~{\rm if}~ i\in I_g.
\end{eqnarray}
Here  
\begin{eqnarray}
\Lambda_i(z)=\frac{T_i(z)}{T_{i+1}(z)}; \quad \tilde{\Lambda}_i(z)=\ln'\Bigg(\frac{T_i(z)T_{i+1}(z)q^{i-1}_+(z)}{q^{i+1}_{+}(z)}\Bigg)\pi_i(z)q^{i+1}_+(z)q^{i-1}_+(z),
\end{eqnarray}
so that $\{\Lambda_i(z)\}$,  $\{T_i(z)\}$ are polynomials and $\pi_i(z)=\prod_k(z-z_k)$, where $z_k$ are distinct roots of  $T_i(z)T_{i+1}(z)$.
The solution of this $qq$-system under nondegeneracy condition is in 1-to-1 correspondence with Bethe ansatz equations for $\mathfrak{gl}(n|m)$ Gaudin model.
The authors of \cite{huang2019diff} produced the $w_i$-composable solutions, which they call populations as in original papers \cite{MV2002}. 
We note, that in this case the extended Weyl group $\tilde{W}$ can be identified with $S_{n+m}$ and in defining representation it can be realized in the standard way using permutation matrices. The architecture of the reproduction procedure is linked to pseudo-differential operator, which for distinguished Dynkin daigram, where one grey root separates Dynkin diagrams for $\mathfrak{sl}(n)$ and  $\mathfrak{sl}(m)$) looks as follows:
\begin{eqnarray}\label{factorber}
R(z)=\prod^n_{i=1}\Bigg(\partial_z-\log'\Bigg[T_i(z)\frac{q^{i-1}_{+}(z)}{q^{i}_{+}(z)}\Bigg]\Bigg)\prod^{m+n}_{n+1}\Bigg(\partial_z+\log'\Bigg[T_i(z)\frac{q^{i-1}_{+}(z)}{q^{i}_{+}(z)}\Bigg]\Bigg)^{-1},
\end{eqnarray}
where $q_+^0(z)=q_+^{m+n}(z)=1$.
\begin{eqnarray}\label{factorberw}
 R^w(z)=
 \prod^{n+m}_{i=1}\Bigg(\partial_z-s_i(w)\log'\Bigg[T^{w}_i(z)\frac{q^{w, i-1}_{+}(z)}{q^{w,i}_{+}(z)}\Bigg]\Bigg)^{s_i(w)}.
\end{eqnarray}
Here $s_i(w)=\pm 1$ and correspond to original permutation $w$ of original $s_i$ in $R(z)$ for distinguished Dynkin diagram. By extra index $w$ we denoted the $qq$-system corresponding to $w$-transformed Dynking diagram. It turns out  that $w$-composability of the related solutions implies identification of $R^w$ and $R$ \cite{huang2019diff}, using the following transformations corresponding to elementary Weyl reflections $w_i$ (this also implies the proper formulas for the $T^w_i(z)$).:
\begin{eqnarray}
&&\Bigg(\partial_z-s_i(w)\log'\Bigg[T^{w}_i(z)\frac{q^{w, i-1}_{+}(z)}{q^{w,i}_{+}(z)}\Bigg]\Bigg)^{s_i(w)}(\partial_z-s_{i+1}(w)\log'\Bigg[T^{w}_{i+1}(z)\frac{q^{w, i}_{+}(z)}{q^{w,i+1}_{+}(z)}\Bigg]\Bigg)^{s_{i+1}(w)}=\nonumber\\
&&\Bigg(\partial_z-s_{i+1}(w)\log'\Bigg[T^{s_iw}_{i}(z)\frac{q^{w, i-1}_{+}(z)}{q^{w,i}_{-}(z)}\Bigg]\Bigg)^{s_{i+1}(w)}
\Bigg(\partial_z-s_i(w)\log'\Bigg[T^{s_iw}_{i+1}(z)\frac{q^{w, i}_{-}(z)}{q^{w,i+1}_{+}(z)}\Bigg]\Bigg)^{s_i(w)}.\nonumber
\end{eqnarray}

This is generalization of original work from \cite{MV2002}, where it was done for the case of $\mathfrak{sl}(n)$ and differential operators. 
At the same time, in that pure even case the shortcut (see \cite{mukhin2009b}, \cite{chertal1}, \cite{chertal2}) between the Gaudin model and opers was achieved by taking a certain analogue of determinant of KZ connection:
\begin{eqnarray}
D^{KZ}_{k,l}=\delta_{k,l}\partial_z-\sum_{\nu=1,\dots, N}\frac{\Phi^\nu_{k,l}}{z-z_\nu}, \quad k,l=1, \dots, n, 
 \end{eqnarray}
where $\Phi^\nu_{k,l}\in U(\mathfrak{gl}(n))^{\otimes N}$ is the operator acting as the $\mathfrak{gl}(n)$ generator $e_{k,l}$  in the $\nu$-th place in that tensor product. This  determinant is understood in a formal way using permutations description. The coefficients in the expansion of the resulting differential operator are the Gaudin Hamiltonians. The factorization of that differential operator of type $R(z)$ corresponds to bringing that operator to the normal form in the sense of \cite{chertal1}, \cite{chertal2},  which coincides with the operator-valued Miura oper connection; in this way the determinant, when applied to a Bethe vector gives exactly the factorization of $R(z)$.

It is natural to expect that the similar formula should exist in the quantum case, since (\ref{factorber}) is very similar to a Berezinian of the corrresponding differential operator. To obtain $R^w$ factorization one can apply the Weyl transformation $w\in \tilde{W}$ represented as a permutation $\tilde{w}$ on $(m|n)$ superspace:
$$
Ber(\tilde{w}[D^{KZ}_{k,l}]\tilde{w}^{-1}).
$$ 
One has to emphasize that $w$ here is not an element of supergroup $SL(m|n)$, but only a formal permutation.  Bringing the $D^{KZ}$-operator to various versions of analogue of normal form for $(m|n)\times (m|n)$ supermatrix should produce an oper connection, which we discussed: after applying formal Berezinian one obtains the $R^w(z)$ operator. 
We mention here a related paper \cite{molevmacmahon} studying center of affine superalgebra $\mathfrak{gl}(n|m)$, where Berezinian formulas naturally emerged; the relation between the center of affine algebras and opers in the pure even case was thorougly investigated (see e.g. \cite{Frenkel_LanglandsLoop} for review).
 
We note, that the generalization of such $qq$-system to the case of 
$\mathfrak{g}=\mathfrak{osp}(m|2n)$ was introduced in \cite{lu2021bethe}, based on  appropriate Bethe equations:

\begin{eqnarray}
&&Wr(q^i_-,q^i_+)(z)=\Lambda_i(z)\prod_i [q^i_+(z)]^{-c_{ij}} ~{\rm if}~  i\in I_w\cup I_b\\\nonumber
&& q^i_-(z)q^i_+(z)=\tilde{\Lambda}_i(z), ~{\rm if}~ i\in I_g.
\end{eqnarray}
Here  
\begin{eqnarray}
\tilde{\Lambda}_i(z)=\log'\Big[\Lambda_i(z)\prod_jq^j_+(z)^{-c_{ij}}\Big]\pi_i(z)\prod_{j, c_{ij}\neq 0}q^j_+(z),
\end{eqnarray} 
where $\{c_{ij}\}$ is the corresponding Cartan matrix, $\Lambda_i(z)\in \mathbb{C}[z]$ for all $i=1, \dots, r$, so that $\pi_i$ is the denominator of the fraction $\log'(\Lambda_i(z))$ of minimal possible degree. 

The construction of the 
$w$-composability (or reproduction procedure) in the language of some pseudo-differential operators similar to $\mathfrak{sl}(n|m)$ case, but much more involved.
 
\subsection{$(G,q)$-opers for supergroups}
There is a generalization of the notion of oper for simple simply-connected Lie groups to $(G,q)$-oper: the difference analogue of G-oper connection, which uses natural multiplicative action of $\mathbb{C}^{\times}_q$ on $\mathbb{P}^1$, i.e. $z\rightarrow qz$\footnote{Such difference connection could be defined for the additive action as well.}. 
 
That$(G,q)$-oper is again a triple 
$(\mathcal{F}_{G}, \mathcal{F}_{B_-}, A)$, where  $\mathcal{F}_{G}, \mathcal{F}_{B_-}$ are the principal $G$-bundle on $\mathbb{P}^1$ and its $B_-$-reduction correspondingly. The $q$-oper connection $A$ is an element of  $Hom(\mathcal{F}_{G}, \mathcal{F^q}_{G})$, where $\mathcal{F^q}_{G}$ is a pull-back bundle with respect to $\mathbb{C}^{\times}_q$-action. Here $A$ belongs to the Coxeter cell: $B_-(z)c(z)B_-(z)$ in the Bruhat decomposition of $G(z)$, where $c(z)$ is a lift of Coxeter element to $G(z)$, which as usual is a  product of the lifts of simple Weyl reflections to $G(z)$.  
The Miura condition adds the requirement that there exists another reduction $\mathcal{F}_{B_+}$ to the Borel subgroup which $A$ preserves.
The $Z$-twist condition means, as in differential case, gauge euivalence of the $q$-connection $A$ to the constant $q$-connection $Z\in H$.

In \cite{Frenkel:2020} we used the notion of $Z$-twisted Miura-Pl\"ucker $(G,q)$-oper to find correspondence with $q$-defomation of the $qq$-system, known as $QQ$-system.  Following that path we constructed the explicit correspondence between $Z$-twisted  $(G,q)$-opers and the XXZ spin chain models XXX spin chains if we use additive action), which are deformations of the Gaudin model,  constructing an example of $q$-Langlands correspondence \cite{KSZ}, \cite{Aganagic:2017smx}. 

Following same principles, one can build a  $q$-oper analogue of ${\bar{\nabla}}$, $\nabla$ for a given simple supergroup and a chosen Dynkin diagram, as an element of:
\begin{eqnarray}\label{sqop}
B_-(z)\Big[\prod_{i\in I_w }{\tilde{w}_{i}(z)}\prod_{j\le k, ~j,k\in I_g\cup I_b}\tilde{w}_{jk}(z)\Big]B_-(z),
\end{eqnarray}
where $B_-$ stands for the Borel subgroup for a chosen Dynkin diagram either for 
supergroup $G$ or the reductive group $\bar{G}$, and 
$\tilde{w}_{i}(z)$, $\tilde{w}_{jk}(z)$ 
stand for the lifts of Weyl reflections corresponding to the roots $\alpha_i$, $\alpha_j+\alpha_k$ correspondingly. 

The explicit expression for Miura$(G,q)$-oper may be quite complicated because of the absence of usual Bruhat decomposition for simple supergroups. At the same time, following  the results of \cite{Frenkel:2020} one can obtain that the corresponding  
$(\bar{G},q)$-oper connection has the following form:
\begin{eqnarray}\label{mqop}
\prod^r_{i=1} [r_i(z)]^{\check{\alpha}_i}\prod_{i\in I_w}e^{R_i(z) e_i}\prod_{j,k\in I_g\cup I_b}e^{R_{ij}(z) [e_j,e_k]},
\end{eqnarray}
where $r_i(z)$, $R_{jk}(z)$ are rational functions.  
Using the expression (\ref{mqop}) one can  relate the notion of $Z$-twisted Miura-Pl\"ucker $(\bar{G},q)$-oper and the $QQ$-systems for simple superalgebras as it is done in the differential case. One may expect that it should lead, using the analogue of the notion of $w$-composable $QQ$-systems to Bethe ansatz equations for $XXX$, $XXZ$ models related to simple superalgebras.  In the case of $\mathfrak{gl}(n|m)$ such examples are also known \cite{huang2019rat}, which are based on the study of rational difference operators: a deformation of the constructions of \cite{huang2019diff}.

\bibliography{cpn1}

\begin{bibdiv}
\begin{biblist}

\bib{Aganagic:2017smx}{article}{
      author={Aganagic, M.},
      author={Frenkel, E.},
      author={Okounkov, A.},
       title={{Quantum $q$-Langlands Correspondence}},
        date={2018},
     journal={Trans. Moscow Math. Soc.},
      volume={79},
       pages={1\ndash 83},
      eprint={1701.03146},
}

\bib{berezin2013introduction}{book}{
      author={Berezin, Feliks~Aleksandrovich},
       title={Introduction to superanalysis},
   publisher={Springer Science \& Business Media},
        date={2013},
      volume={9},
}

\bib{Beilinson:2005}{article}{
      author={Beilinson, A.},
      author={Drinfeld, V.},
       title={{Opers}},
        date={2005},
      eprint={math/0501398},
         url={https://arxiv.org/abs/math/0501398},
}

\bib{bergvelt1999supercurves}{article}{
      author={Bergvelt, M.J.},
      author={Rabin, J.M.},
       title={Supercurves, their jacobians, and super kp equations},
        date={1999},
     journal={Duke Math. J.},
      volume={98},
       pages={1\ndash 57},
}

\bib{baranov1985multiloop}{article}{
      author={Baranov, M.A.},
      author={Shvarts, A.S.},
       title={Multiloop contribution to string theory},
        date={1985},
     journal={JETP Lett.},
      volume={42},
       pages={419\ndash 421},
}

\bib{Brinson:2021ww}{article}{
      author={Brinson, Ty~J.},
      author={Sage, Daniel~S.},
      author={Zeitlin, Anton~M.},
       title={{Opers on the projective line, Wronskian relations, and the Bethe
  Ansatz}},
      eprint={2112.02711},
         url={https://arxiv.org/pdf/2112.02711.pdf},
}

\bib{chertal1}{article}{
      author={Chervov, A.},
      author={Talalaev, D.},
       title={{Quantum spectral curves, quantum integrable systems and the
  geometric Langlands correspondence}},
      eprint={hep-th/0604128},
}

\bib{chertal2}{article}{
      author={Chervov, A.},
      author={Talalaev, D.},
       title={{KZ equation, G-opers, quantum Drinfeld--Sokolov reduction, and
  quantum Cayley--Hamilton identity}},
        date={2008},
     journal={Zap. Nauchn. Sem. POMI},
      volume={360},
       pages={246\ndash 259},
      eprint={hep-th/0607250},
}

\bib{cheng2012dualities}{book}{
      author={Cheng, Shun-Jen},
      author={Wang, Weiqiang},
       title={{Dualities and representations of Lie superalgebras}},
   publisher={AMS},
        date={2012},
}

\bib{delduc1998supersymmetric}{article}{
      author={Delduc, F.},
      author={Gallot, Laurent},
       title={{Supersymmetric Drinfeld--Sokolov reduction}},
        date={1998},
     journal={J. Math. Phys.},
      volume={39},
      number={9},
       pages={4729\ndash 4745},
      eprint={solv-int/9802013},
}

\bib{Frenkel:1995zp}{article}{
      author={Frenkel, E.},
       title={{Affine algebras, Langlands duality and Bethe ansatz}},
        date={1995},
     journal={{Proceedings of the International Congress of Mathematical
  Physics, Paris, 1994, ed. D. Iagolnitzer, pp. 606--642, International
  Press.}},
      eprint={q-alg/9506003},
}

\bib{Frenkel:2003qx}{article}{
      author={Frenkel, E.},
       title={{Opers on the Projective Line, flag Manifolds and Bethe Ansatz}},
        date={2003},
     journal={Mosc. Math. J.},
      volume={4},
       pages={655\ndash 705},
      eprint={math/0308269},
}

\bib{Frenkel:2004qy}{article}{
      author={Frenkel, E.},
       title={{Gaudin Model and Opers}},
        date={2004},
     journal={{Infinite Dimensional Algebras and Quantum Integrable Systems,
  eds. P. Kulish, e.a., Progress in Math.}},
      volume={237},
       pages={1\ndash 60},
      eprint={math/0407524},
}

\bib{Frenkel_LanglandsLoop}{book}{
      author={Frenkel, Edward},
       title={Langlands correspondence for loop groups},
   publisher={Cambridge University Press},
        date={2007},
}

\bib{Feigin:1994in}{article}{
      author={Feigin, B.},
      author={Frenkel, E.},
      author={Reshetikhin, N.},
       title={{Gaudin Model, Bethe Ansatz and Critical Level}},
        date={1994},
     journal={Commun. Math. Phys.},
      volume={166},
       pages={27\ndash 62},
      eprint={hep-th/9402022},
}

\bib{Rybnikov:2010}{article}{
      author={Feigin, B.},
      author={Frenkel, E.},
      author={Rybnikov, L.},
       title={{Opers with irregular singularity and spectra of the shift of
  argument subalgebra}},
        date={2010},
     journal={Duke Math. J.},
      volume={155},
       pages={337\ndash 363},
      eprint={0712.1183},
}

\bib{friedan1986notes}{article}{
      author={Friedan, D.},
       title={Notes on string theory and two-dimensional conformal field
  theory},
        date={1986},
     journal={{In: The proceedings of the workshop on unified string theories.
  Gross, D., Green, M. (eds.), Singapore: World Press}},
}

\bib{Feigin:2006xs}{article}{
      author={Feigin, B.},
      author={Frenkel, E.},
      author={Toledano~Laredo, V.},
       title={{Gaudin models with irregular singularities}},
        date={2010},
     journal={Adv. Math.},
      volume={223},
       pages={873\ndash 948},
      eprint={math/0612798},
}

\bib{Frenkel:2020}{article}{
      author={Frenkel, Edward},
      author={Koroteev, Peter},
      author={Sage, Daniel~S.},
      author={Zeitlin, Anton~M.},
       title={{q-Opers, QQ-Systems, and Bethe Ansatz}},
     journal={J. Eur. Math. Soc.},
       pages={to appear},
      eprint={2002.07344},
}

\bib{frappat2000dictionary}{book}{
      author={Frappat, Luc},
      author={Sciarrino, Antonino},
      author={Sorba, Paul},
       title={Dictionary on lie algebras and superalgebras},
   publisher={Academic Press San Diego, CA},
        date={2000},
      volume={10},
}

\bib{gualzetti1993quantum}{article}{
      author={Gualzetti, Alessandro},
      author={Penati, Silvia},
      author={Zanon, Daniela},
       title={{Quantum conserved currents in supersymmetric Toda theories}},
        date={1993},
     journal={Nucl. Phys. B},
      volume={398},
      number={3},
       pages={622\ndash 658},
}

\bib{huang2019rat}{article}{
      author={Huang, Chenliang},
      author={Lu, Kang},
      author={Mukhin, Evgeny},
       title={{Solutions of XXX Bethe ansatz equation and rational difference
  operators}},
        date={2019},
     journal={J. Phys. A},
      volume={52},
      number={37},
       pages={375204},
      eprint={1811.11225},
}

\bib{huang2019diff}{article}{
      author={Huang, Chenliang},
      author={Mukhin, Evgeny},
      author={Vicedo, Beno{\^\i}t},
      author={Young, Charles},
       title={{The solutions of $\mathfrak{gl}_{M|N}$ Bethe ansatz equation and
  rational pseudodifferential operators}},
        date={2019},
     journal={Sel. Math. New Ser.},
      volume={25},
       pages={1\ndash 34},
      eprint={1809.01279},
}

\bib{inami1991lie}{article}{
      author={Inami, Takeo},
      author={Kanno, Hiroaki},
       title={{Lie superalgebraic approach to super Toda lattice and
  generalized super KdV equations}},
        date={1991},
     journal={Commun. Math. Phys.},
      volume={136},
      number={3},
       pages={519\ndash 542},
}

\bib{kac2006representations}{inproceedings}{
      author={Kac, V.G.},
       title={{Representations of classical Lie superalgebras}},
organization={Springer},
   booktitle={{Differential Geometrical Methods in Mathematical Physics II:
  Proceedings, University of Bonn, July 13--16, 1977}},
       pages={597\ndash 626},
}

\bib{kac1977lie}{article}{
      author={Kac, Victor~G.},
       title={{Lie superalgebras}},
        date={1977},
     journal={Adv. Math.},
      volume={26},
      number={1},
       pages={8\ndash 96},
}

\bib{kulish2001bethe}{article}{
      author={Kulish, P.P.},
      author={Manojlovi{\'c}, N.},
       title={{Bethe vectors of the $osp(1|2)$ Gaudin model}},
        date={2001},
     journal={Lett. Math. Phys.},
      volume={55},
       pages={77\ndash 95},
}

\bib{kulish2001creation}{article}{
      author={Kulish, P.P.},
      author={Manojlovi{\'c}, N.},
       title={{Creation operators and Bethe vectors of the $osp (1| 2)$ Gaudin
  model}},
        date={2001},
     journal={J. Math. Phys.},
      volume={42},
      number={10},
       pages={4757\ndash 4778},
}

\bib{KSZ}{article}{
      author={Koroteev, P.},
      author={Sage, D.~S.},
      author={Zeitlin, A.M.},
       title={{(SL(N),q)-Opers, the q-Langlands Correspondence, and
  Quantum/Classical Duality}},
        date={2021},
     journal={Commun. Math. Phys.},
      volume={381},
      number={2},
       pages={641\ndash 672},
      eprint={1811.09937},
         url={https://arxiv.org/abs/1508.02690},
}

\bib{KZminors}{article}{
      author={Koroteev, P.},
      author={Zeitlin, A.M.},
       title={{q-Opers, QQ-systems, and Bethe Ansatz II: Generalized Minors}},
        date={2023},
     journal={J. Reine Angew. Math.},
      volume={795},
       pages={271\ndash 296},
      eprint={2108.04184},
}

\bib{Koroteev:2021a}{article}{
      author={Koroteev, Peter},
      author={Zeitlin, Anton~M.},
       title={{3d Mirror Symmetry for Instanton Moduli Spaces}},
      eprint={2105.00588},
         url={https://arxiv.org/pdf/2105.00588.pdf},
}

\bib{Koroteev2022zoo}{article}{
      author={Koroteev, Peter},
      author={Zeitlin, Anton~M.},
       title={{The Zoo of Opers and Dualities}},
      eprint={2208.08031},
}

\bib{Koroteev:2020tv}{article}{
      author={Koroteev, Peter},
      author={Zeitlin, Anton~M.},
       title={Toroidal q-opers},
        date={2023},
     journal={J. Inst. Math. Jussieu},
      volume={22},
       pages={581\ndash 642},
      eprint={2007.11786},
}

\bib{leitessup}{book}{
      author={Leites, D.A.},
       title={{Theory of supermanifolds}},
   publisher={KF Akad. Nauk SSSR, Petrozavodsk},
        date={1983},
}

\bib{lu2021bethe}{article}{
      author={Lu, Kang},
      author={Mukhin, Evgeny},
       title={{Bethe ansatz equations for orthosymplectic Lie superalgebras and
  self-dual superspaces}},
organization={Springer},
        date={2021},
     journal={Ann. Henri Poincar{\'e}},
      volume={22},
      number={12},
       pages={4087\ndash 4130},
      eprint={2103.16729},
}

\bib{leites1985embeddings}{inproceedings}{
      author={Leites, D.A.},
      author={Saveliev, M.V.},
      author={Serganova, V.V.},
       title={{Embedding of $\mathfrak{osp}(N|2)$ and the associated nonlinear
  supersymmetric equations}},
        date={1986},
   booktitle={{Group theoretical methods in physics, Vol. I (Yurmala, 1985) }},
   publisher={VNU Sci. Press, Utrecht},
       pages={255\ndash 297},
}

\bib{manin2013gauge}{book}{
      author={Manin, Yuri~I},
       title={Gauge field theory and complex geometry},
   publisher={Springer Science \& Business Media},
        date={2013},
      volume={289},
}

\bib{molevmacmahon}{article}{
      author={Molev, Alexander~Ivanovich},
      author={Ragoucy, Eric},
       title={{The MacMahon Master Theorem for right quantum superalgebras and
  higher Sugawara operators for $\mathfrak{gl}(m|n)$}},
        date={2014},
     journal={Mosc.Math. J.},
      volume={14},
       pages={83\ndash 119},
      eprint={0911.3447},
}

\bib{mukhin2009b}{article}{
      author={Mukhin, Evgeny},
      author={Tarasov, Vitaly},
      author={Varchenko, Alexander},
       title={{The B. and M. Shapiro conjecture in real algebraic geometry and
  the Bethe ansatz}},
        date={2009},
     journal={Annals of mathematics},
      volume={170},
      number={2},
       pages={863\ndash 881},
      eprint={math/0512299},
}

\bib{mukhvarmiura}{article}{
      author={Mukhin, E.},
      author={Varchenko, A.},
       title={{Miura Opers and Critical Points of Master Functions}},
     journal={Cent. Eur. J. Math.},
      volume={3},
      eprint={math/0312406},
         url={https://arxiv.org/abs/math/0312406},
}

\bib{MV2002}{article}{
      author={Mukhin, E.},
      author={Varchenko, A.},
       title={Critical points of master functions and flag varieties},
        date={2004},
     journal={Communications in Contemporary Mathematics},
      volume={6},
      number={1},
       pages={111\ndash 163},
      eprint={math/0209017},
         url={https://arxiv.org/abs/math/0209017},
}

\bib{mukhvar2008}{article}{
      author={Mukhin, E.},
      author={Varchenko, A.},
       title={{Quasi-polynomials and the Bethe Ansatz }},
        date={2008},
     journal={Geom. Top. Mon.},
      volume={13},
      eprint={math/0604048},
}

\bib{mukhin2015gaudin}{article}{
      author={Mukhin, Evgeny},
      author={Vicedo, Beno{\^\i}t},
      author={Young, Charles},
       title={{Gaudin models for $\mathfrak{gl}(m|n)$}},
        date={2015},
     journal={J. Math. Phys.},
      volume={56},
      number={5},
       pages={051704},
}

\bib{olshanetsky1983supersymmetric}{article}{
      author={Olshanetsky, M.A.},
       title={{Supersymmetric two-dimensional Toda lattice}},
        date={1983},
     journal={Commun. Math. Phys.},
      volume={88},
       pages={63\ndash 76},
}

\bib{penkov1990borel}{article}{
      author={Penkov, I.B.},
       title={{Borel-Weil-Bott theory for classical Lie superalgebras}},
        date={1990},
     journal={J. Sov. Math.},
      volume={51},
       pages={2108\ndash 2140},
}

\bib{RT}{article}{
      author={Rakowski, Mark},
      author={Thompson, George},
       title={Connections on vector bundles over super riemann surfaces},
        date={1989},
     journal={Phys. Lett. B},
      volume={220},
      number={4},
       pages={557\ndash 561},
}

\bib{manvor2}{article}{
      author={Voronov, A.A.},
      author={Manin, Yu.~I.},
       title={Supercell partitions of flag superspaces},
        date={1990},
     journal={J. Sov. Math.},
      volume={51},
       pages={2083\ndash 2108},
}

\bib{manvor1}{article}{
      author={Voronov, Alexander~A.},
      author={Manin, Yu.~I.},
       title={Schubert supercells},
        date={1984},
     journal={Funct. Anal. Appl.},
      volume={18},
      number={4},
       pages={329\ndash 330},
}

\bib{witten2012notes}{article}{
      author={Witten, Edward},
       title={{Notes on super Riemann surfaces and their moduli}},
      eprint={1209.2459},
}

\bib{zeitlin2022wronskians}{article}{
      author={Zeitlin, A.M.},
       title={{On Wronskians and $qq$-systems}},
     journal={to appear in Contemp. Math.},
      eprint={2208.08018},
}

\bib{Zeitlin:2013iya}{article}{
      author={Zeitlin, A.M.},
       title={{Superopers on Supercurves}},
        date={2015},
     journal={Lett. Math. Phys.},
      volume={105},
      number={2},
       pages={149\ndash 167},
      eprint={1311.5997},
}

\end{biblist}
\end{bibdiv}

\end{document}